\documentclass{ws-procs9x6}

\usepackage{graphicx}

\def\R{\mathbb{R}}

\def\Image{\mathop\mathrm{Image}}

\def\<{\langle}
\def\>{\rangle}

\begin{document}

\title{NUMERICAL ASPECTS OF EVOLUTION OF PLANE CURVES SATISFYING THE FOURTH ORDER GEOMETRIC EQUATION}

\author{K. MIKULA$^{*}$ and D. \v SEV\v COVI\v C$^{**}$ }

\address{
$^{*}$
Department of Mathematics and Descriptive Geometry, Slovak University of
Technology, 813 68 Bratislava, Slovak Republic,
\\
E-mail: mikula@math.sk
\\[5pt]
$^{**}$
Faculty of Mathematics, Physics and Informatics, 
Comenius University, 
\\
842 48 Bratislava, Slovak Republic,
\\
E-mail: sevcovic@fmph.uniba.sk
}

\begin{abstract}
In this review paper we present a stable Lagrangian numerical method 
for computing plane curves evolution driven by  the fourth order 
geometric equation.  The numerical scheme and 
computational examples are presented.
\end{abstract}


\bodymatter

\section{ Introduction}
\label{sec-1}
The main purpose of this contribution is to suggest a method for computing evolution of 
closed smooth plane curves  driven by the normal velocity $\beta$ depending
on the intrinsic Laplacian  of the curvature $k$ and curvature itself:
\begin{equation}
\beta = - \partial^2_s k  + b(k)
\label{geomrov}
\end{equation}
where $b(k)$ is a $C^2$ function of the curvature $k$, $b(0)=0$.
Numerical aspects of evolution of plane curves satisfying (\ref{geomrov}) have been studied in \cite{MS_ALG} for the case of the surface diffusion flow with no lower order terms, i.e.  $\beta = - \partial^2_s k$, and in \cite{BMOS} for the case of the so-called Willmore flow for which the normal velocity is given by $\beta = - \partial^2_s k -\frac12 k^3$.
Recall that the latter case corresponds to the motion of elastic curves, e.g. the model 
of Euler-Bernoulli elastic rod, which  is an important problem in structural mechanics.
The elastic curve evolution and surface diffusion can be found 
in many practical applications as sintering (in brick production),
formation of rock strata from sandy sediments, metal thin film growth etc.
(see e.g. \cite{Se2}). 

\section{Governing equations}

We follow the so-called direct (or Lagrangian)  approach. We represent a solution of (\ref{geomrov}) by 
the position vector $x$ satisfying  the geometric equation $\partial_t x = \beta \vec N + \alpha \vec T$ where $\vec N, 
\vec T$ are the unit inward normal and  tangent vectors. 
An immersed regular plane  curve $\Gamma$ can be parameterized by a smooth function 
$x:S^1\to \R^2$, i.e. $\Gamma=\{ x(u), u\in S^1 \}$, for which  
$g=|\partial_u x| >0$. Taking into account the periodic boundary conditions at 
$u=0,1$ we shall hereafter identify $S^1$ with the interval $[0,1]$. The unit 
arc-length parameterization will be denoted by 
$s$, so $\mbox{d} s = g\, \mbox{d} u$. 
The tangent vector $\vec T$ and the signed curvature 
$k$ of $\Gamma$ satisfy 
$\vec T = \partial_s x$,
$k = \mbox{det}(\partial_s x, \partial^2_s x)$.
Moreover, we choose the unit inward normal vector $\vec N$  such that 
$\mbox{det}(\vec T,\vec N) =1$. 
Let a regular smooth initial curve $\Gamma_0=\Image(x_0)$ be given. It turns out that a family 
of plane curves $\Gamma_t = \Image(x(., t)), t\in [0,T)$, satisfying 
(\ref{geomrov}) can be represented  by a solution 
to the following system of PDEs:
\begin{eqnarray}
&&\partial_t k  =  - \partial^4_s k    +\partial^2_s b(k)
+ \partial_s(\alpha k) + k(k\beta - \partial_s \alpha),
\label{rov2} \\
&&\partial_t \eta = - k\beta + \partial_s \alpha,\ \ \ g=\exp(\eta),
\label{rov1}\\
&&\partial_t x = -\partial_s^4 x   - \frac{k^3-b(k)}{k}\partial^2_s x  
+ (\alpha- \frac32 \partial_s (k^2) ) \partial_s x\,,
\label{rov4}
\end{eqnarray}
subject to initial conditions
$k(.,0)=k_0\,,\   g(.,0)=g_0\,, \ x(.,0)=x_0(.)$ 
and periodic boundary conditions at $u=0,1$ (cf. \cite{MS2,MS_CVS}).

\section{Approximation scheme and numerical experiments}

\begin{figure}
\vglue -0.5truecm
\centerline{
\includegraphics[width=3truecm]{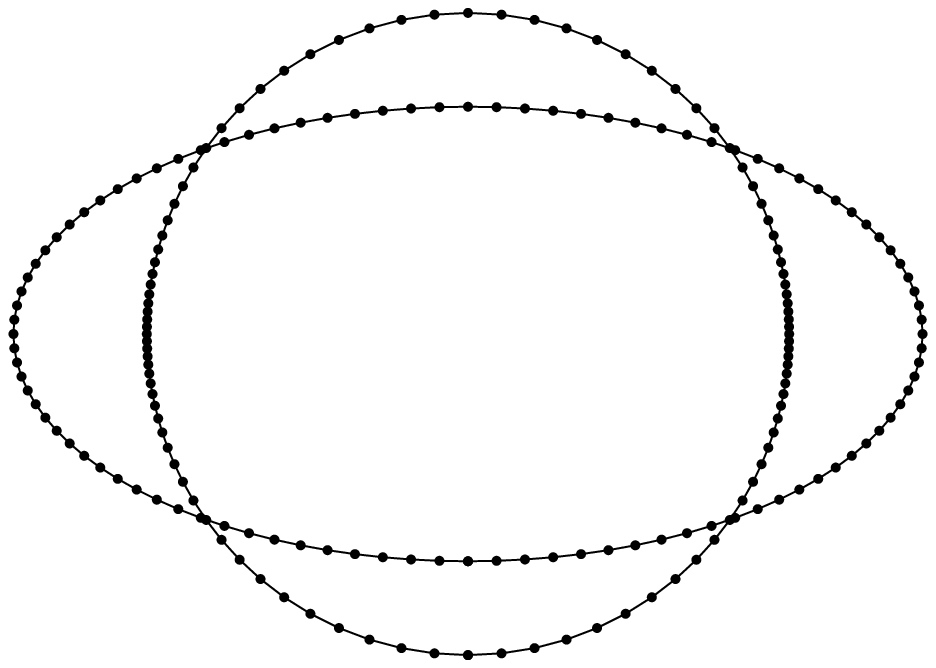}
\hglue 0.3cm
\includegraphics[width=3truecm]{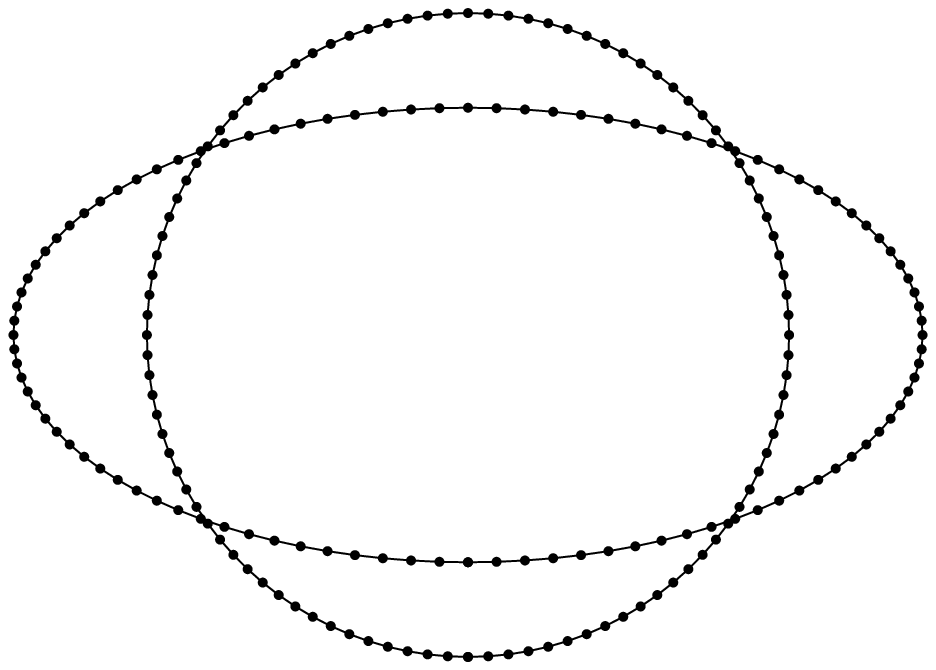}
\hglue 0.3cm
\includegraphics[width=2.8cm]{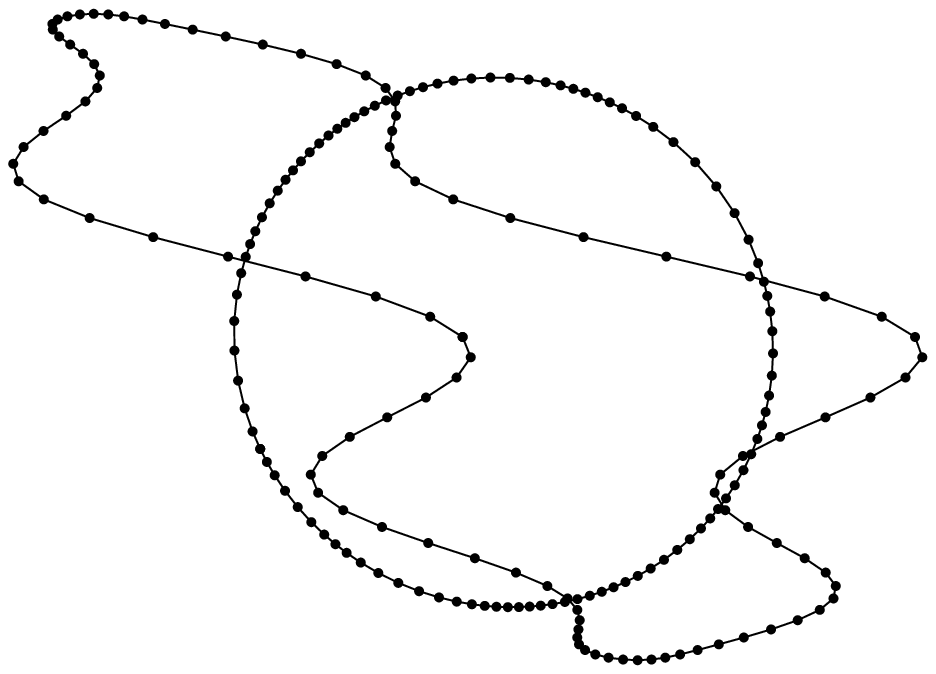}
}
\centerline{a) \hglue3truecm b) \hglue 3truecm c)}
\caption{
Surface diffusion of an initial ellipse and its circular asymptotic
shapes without a) and with b) tangential redistribution. 
Evolution of the highly nonconvex initial curve and its  
asymptotic circular shape at time $t=0.170$ c).}
\label{fig1}
\end{figure}

\begin{figure}
\vglue-1truecm
\center{
\includegraphics[width=2.5cm]{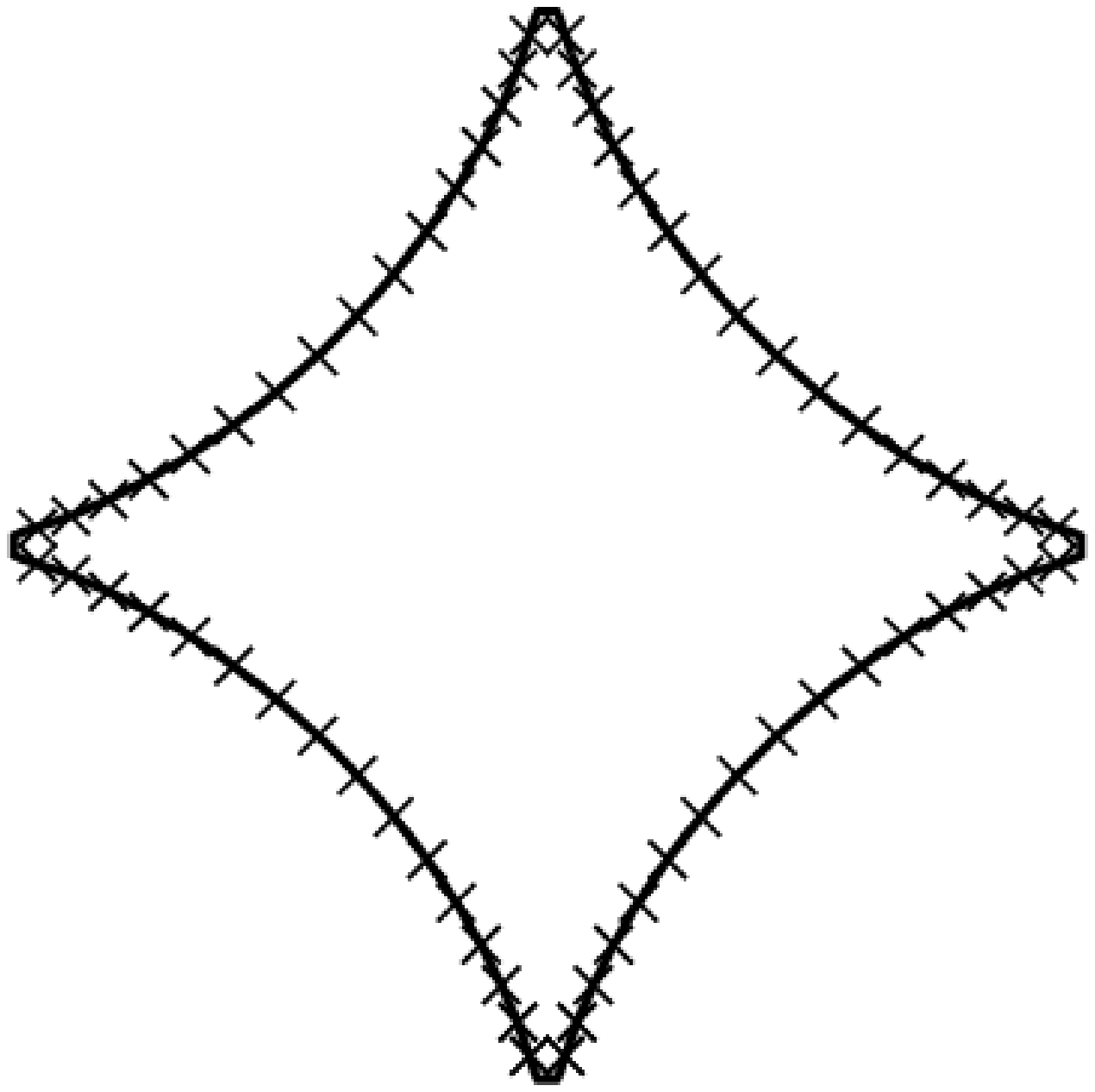}
\hglue 0.5cm
\includegraphics[width=2.5cm]{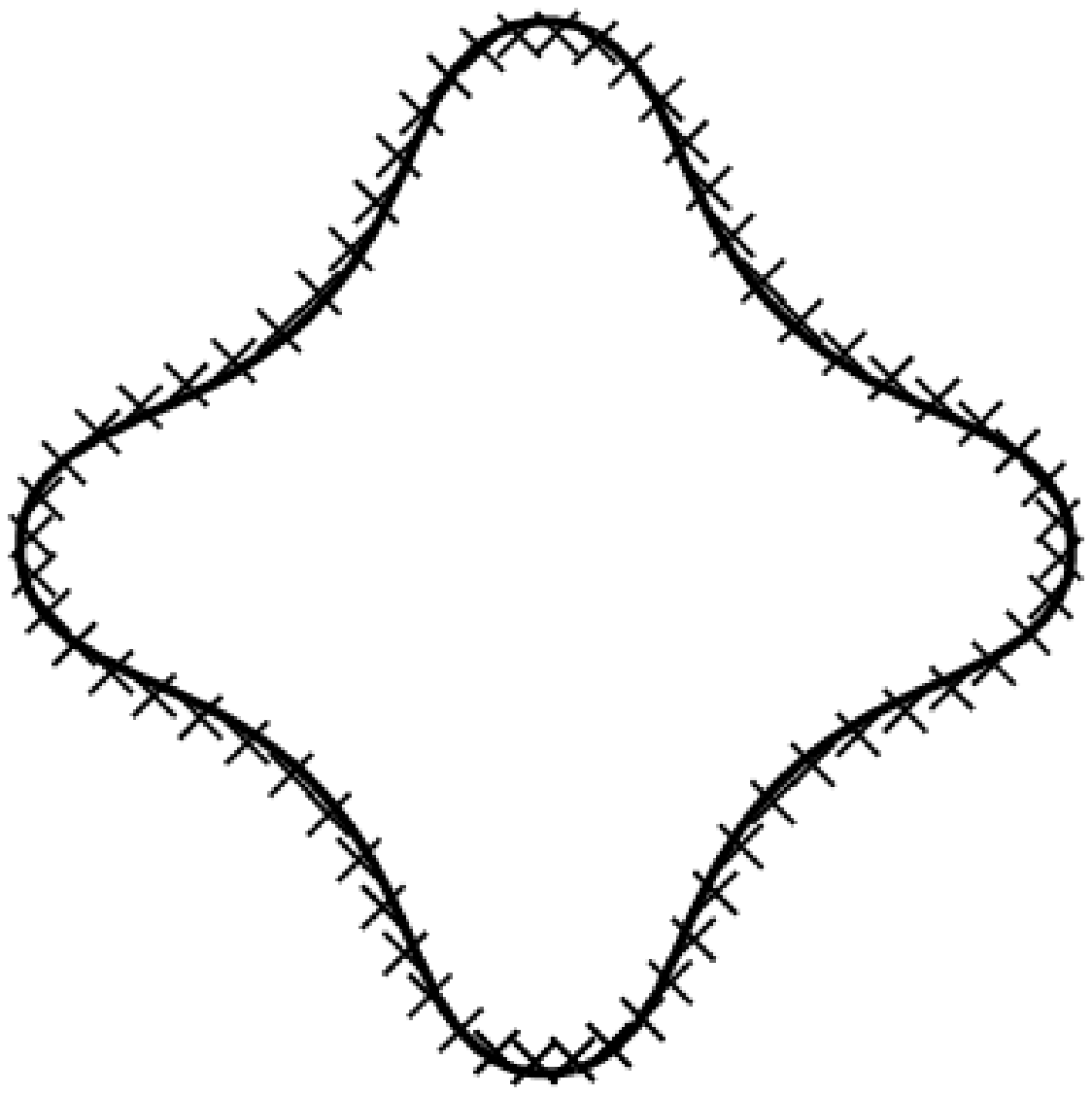}
\hglue 0.5cm
\includegraphics[width=2.5cm]{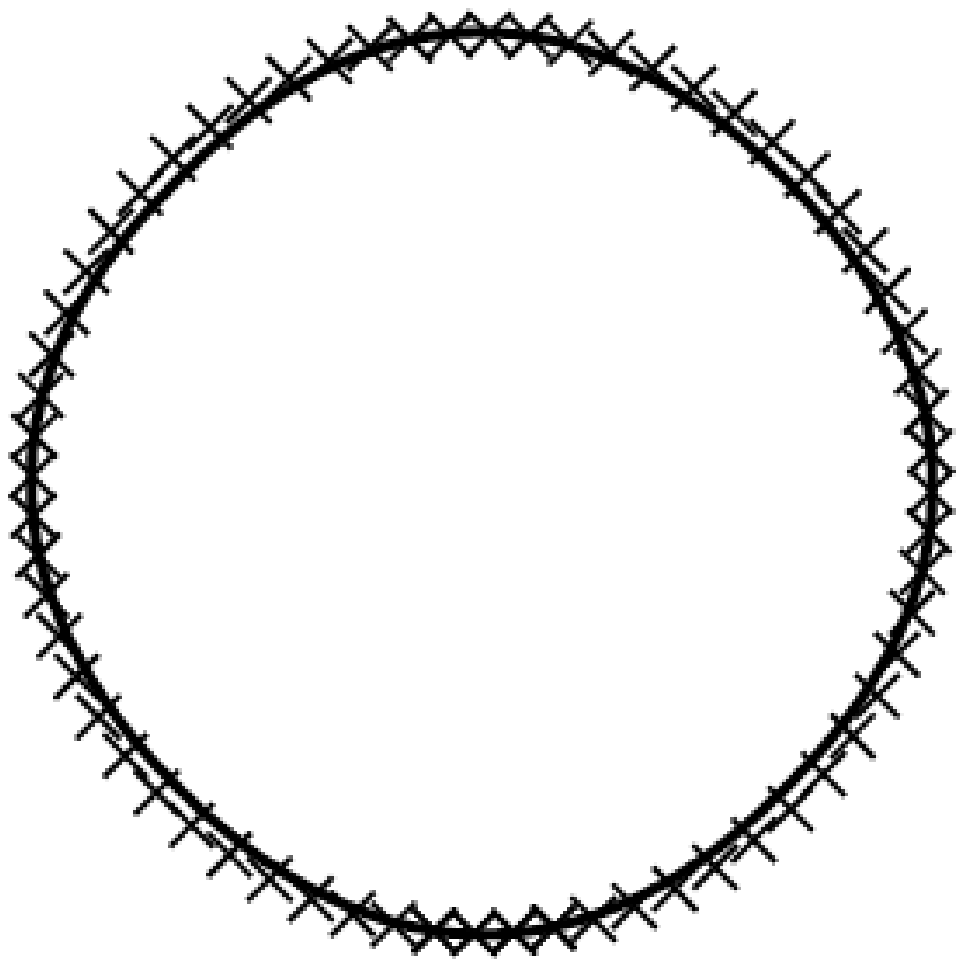}
}
\caption{
An asteroid as an initial condition and its evolution by the Willmore flow
at times  
$t=0,\, 0.0005, \, 0.005$.}
\label{astroid-test}
\end{figure}

The presence of a tangential velocity $\alpha$  in the position vector equation has no 
impact on the shape of evolving curves. As it was shown e.g. in \cite{De,D2,Hou1, K2,MS2, MS_CVS} 
for general curvature driven motions (nonlinear, anisotropic, with external forces) incorporation of a suitable tangential velocity into governing equations stabilizes numerical computations significantly.
It prevents the direct Lagrangian algorithm from its main drawbacks -- the merging of numerical grid points and their order exchange. It also allows for larger time steps without loosing stability. In our numerical solution we consider tangential velocity given by a  nonlocal tangential redistributions discussed in a detail in \cite{Hou1,K2,MS2,MS_CVS}. It follows from \cite{MS2,MS_CVS} that the redistribution functional $\alpha$ satisfying
\begin{equation}
\partial_s \alpha = k\beta -\langle k\beta \rangle_\Gamma
+ \omega \left(L/g -1\right)\,,
\label{alpha-general}
\end{equation}
with a constant $\omega>0$, 
is capable of asymptotic uniform redistribution of grid points along 
the evolved curve. 
In our computational method a numerically evolved curve is represented by discrete plane points 
$x_i^j$ where the index $i=1,...,n,$ denotes space 
discretization and the index $j=0,...,m,$ denotes a discrete time stepping. 
The linear approximation of an evolving curve in the $j$-th discrete time step 
is thus given by a polygon with vertices $x_i^j, i=1,...,n$. Due to
periodicity conditions we shall also use additional values $x_{-1}^j=x_{n-1}^j$, 
$x_0^j=x_n^j$, $x_{n+1}^j=x_1^j$, $x_{n+2}^j=x_2^j$. 
If we take a uniform division of the time interval
$[0,T]$ with a time step $\tau=\frac{T}{m}$ and a uniform division of the fixed
parameterization interval $[0,1]$ with a step $h=\frac1n$, a point $x_i^j$
corresponds to $x(ih,j\tau)$. 
The systems of difference equations corresponding to 
(\ref{rov2})--(\ref{rov4}) and (\ref{alpha-general}) will be given 
for discrete quantities $\alpha_i^j$, $\eta_i^j$, $r_i^j$, $k_i^j$, 
$x_i^j$, $i=1,...,n,\ j=1,...,m$, representing approximations of the unknowns $\alpha$, 
$\eta$, $g h$, $k$,  and $x$, respectively. Here $\alpha_i^j$ represents 
tangential velocity of a flowing node $x_i^j$, and $\eta_i^j$,
$r_i^j\approx |x_i^j-x_{i-1}^j|$, $k_i^j$, $\nu_i^j$ 
represent piecewise constant approximations of 
the corresponding quantities in the so-called 
{\it flowing finite volume} $\left[x_{i-1}^j,x_i^j\right]$. 
We shall use the corresponding flowing dual volumes
$\left[\tilde x_{i-1}^j,\tilde x_i^j\right]$, where 
$\tilde x_i^j= \frac{x_{i-1}^j+x_{i}^j}{2}$, with
approximate lengths
$q_i^j\approx |\tilde x_i^j-\tilde x_{i-1}^j|$.
At the $j$-th discrete time step, we first find discrete values of the 
tangential velocity $\alpha_i^j$ by discretizing equation (\ref{alpha-general}). 
Then the values of redistribution parameter $\eta_i^j$ are computed
and utilized for updating discrete local lengths $r_i^j$ by discretizing equations 
(\ref{rov1}). Using already computed 
local lengths, the intrinsic derivatives are approximated in (\ref{rov2}),
 and (\ref{rov4}), and pentadiagonal systems with 
periodic boundary conditions are constructed and solved for discrete 
curvatures $k_i^j$ and position vectors $x_i^j$. 
In the sequel, we present in a more detail our discretization.
Using $r_i^{j-1}$ as an approximation of the 
length of the flowing finite volume 
$\left[x_{i-1}^{j-1},x_i^{j-1}\right]$ at the previous $j$th time step
we construct difference approximation of the
intrinsic derivative $\partial_s\alpha\approx \frac
{\alpha_i^j-\alpha_{i-1}^j}{r_i^{j-1}}$ and by taking all further quantities in 
(\ref{alpha-general}) from the previous time step. We obtain
the following expression for {\it discrete values of the tangential
velocity:}
$\alpha_i^j = \alpha_{i-1}^j +
 r_i^{j-1}k_i^{j-1}\beta_i^{j-1} - r_i^{j-1} B^{j-1}  
+ \omega(M^{j-1}- r_i^{j-1}) $ 
where
$\beta_i^{j-1}= - \frac{1}{r_i^{j-1}}
\left(\frac{k_{i+1}^{j-1}-k_i^{j-1}}{q_i^{j-1}}-
\frac{k_{i}^{j-1}-k_{i-1}^{j-1}}{q_{i-1}^{j-1}}\right)+b(k_i^{j-1})$, 
$q_i^{j-1}=\frac{r_{i}^{j-1}+r_{i+1}^{j-1}}{2}$, $M^{j-1}=\frac{1}{n} L^{j-1}$,  
$L^{j-1}=\sum\limits_{l=1}^n r_l^{j-1}$, $B^{j-1}=\frac{1}{L^{j-1}}
\sum\limits_{l=1}^n  r_l^{j-1}k_l^{j-1}\beta_l^{j-1}$ and $\alpha_0^j=0$, 
i.e. the point $x_0^j$ is moved in the normal direction only. 
Inserting (\ref{alpha-general}) in (\ref{rov1}) and using 
a similar strategy  
give us: $r_i^{j-1}\frac{\eta_i^j-\eta_i^{j-1}}{\tau} =
- r_i^{j-1} B^{j-1}+\omega(M^{j-1}- r_i^{j-1})$,
for $i=1,...,n, \eta_0^j=\eta_n^j,\ \ \eta_{n+1}^j=\eta_1^j$. 
Next we {\it update local lengths} by the rule: $r_i^j=\exp(\eta_i^j),  r_{-1}^j=r_{n-1}^j,
\ \ r_0^j=r_n^j,\ \ r_{n+1}^j=r_1^j, \ \ r_{n+2}^j=r_2^j.$
Subsequently,  new local lengths are used for approximation of intrinsic 
derivatives in (\ref{rov2})  and (\ref{rov4}). 
First, we derive a discrete analogy of the curvature equation (\ref{rov2}). 
We have to approximate the 4-th order derivative of curvature inside the
flowing finite volume $\left[x_{i-1},x_i\right]$, $i=1,...,n$. 
For that goal we take the following finite difference approximation:
$\partial_s^4 k(\tilde x_i)  \approx 
\frac 1{r_i q_i r_{i+1} q_{i+1}} k_{i+2} -
\bigl(\frac 1{r_i q_i r_{i+1} q_{i+1}} + \frac 1{r_i q_i^2 r_{i+1}} +
\frac 1{r_i^2 q_i^2}+\frac 1{r_i^2 q_i q_{i-1}}\bigr) k_{i+1} +
\bigl(\frac 1{r_i q_i^2 r_{i+1}} +
\frac 1{r_i^2 q_i^2}+\frac 2{r_i^2 q_i q_{i-1}}
+\frac 1{r_i^2 q_{i-1}^2}+\frac 1{r_i q_{i-1}^2 r_{i-1}}\bigr) k_i +
\bigl(\frac 1{r_i^2 q_i q_{i-1}}
+\frac 1{r_i^2 q_{i-1}^2}+\frac 1{r_i q_{i-1}^2 r_{i-1}}
+\frac 1{r_i q_{i-1} r_{i-1} q_{i-2}}\bigr) k_{i-1} +
\frac 1{r_i q_{i-1} r_{i-1} q_{i-2}} k_{i-2}$. 
Approximating first and second order terms in (\ref{rov2}) 
by central differences
and taking semi-implicit time stepping we obtain following {\it pentadiagonal 
system} with periodic boundary conditions for {\it new 
discrete values of curvature}:
\[
a_{i}^j k_{i-2}^j+b_{i}^j k_{i-1}^j+c_{i}^j k_{i}^j+d_i^j k_{i+1}^j
+e_i^j k_{i+2}^j = f_i^j
\]
subject to periodic b.c. 
$k_{-1}^j=k_{n-1}^j, k_0^j=k_n^j, k_{n+1}^j=k_1^j,  k_{n+2}^j=k_2^j,$
where
\begin{eqnarray*}
&&a_i^j = \frac 1{q_{i-1}^j r_{i-1}^j q_{i-2}^j},\ \ 
e_i^j = \frac 1{q_i^j r_{i+1}^j q_{i+1}^j},
\\ 
&&f_i^j = \frac {r_i^j}{\tau}k_i^{j-1}
+\frac{b(k_{i+1}^{j-1})-b(k_{i}^{j-1})}{q_i^j}
-
\frac{b(k_{i}^{j-1}) -b(k_{i-1}^{j-1})}{q_{i-1}^j}
,\nonumber\\
&&b_i^j = -\left(\frac 1{r_i^j q_i^j q_{i-1}^j}
+\frac 1{r_i^j (q_{i-1}^j)^2}+\frac 1{(q_{i-1}^j)^2 r_{i-1}^j}
+\frac 1{q_{i-1}^j r_{i-1}^j q_{i-2}^j}\right)+\frac{\alpha_{i-1}^j}{2}
\nonumber\\
&&d_i^j = -\left(\frac 1{q_i^j r_{i+1}^j q_{i+1}^j} + 
\frac 1{(q_i^j)^2 r_{i+1}^j} +
\frac 1{r_i^j (q_i^j)^2}+\frac 1{r_i^j q_i^j q_{i-1}^j}\right)
-\frac{\alpha_{i}^j}{2}
\nonumber\\
&&c_i^j = \frac 1{(q_i^j)^2 r_{i+1}^j} +
\frac 1{r_i^j (q_i^j)^2}+\frac 2{r_i^j q_i^j q_{i-1}^j}
+\frac 1{r_i^j (q_{i-1}^j)^2}+\frac 1{(q_{i-1}^j)^2 r_{i-1}^j}+
\nonumber\\
&&\hskip 1truecm \frac{r_i^j}{\tau}
- r_i^{j-1} k_i^{j-1}\beta_i^{j-1}+\frac{\alpha_i^j}{2}-\frac{\alpha_{i-1}^j}{2}.
\nonumber
\end{eqnarray*}
In order to construct discretization of equation (\ref{rov4})
we approximate the intrinsic derivatives  in a dual volume
$\left[\tilde x_{i-1},\tilde x_i\right]$. For approximation of the fourth order
intrinsic derivative of the position vector we take similar approach as
above for curvature, but in the middle point $x_i$ of the dual volume.
In such a way and using the semi-implicit approach,
we end up with {\it two tridiagonal systems for updating 
the discrete position vector:}
\[
{\cal A}_{i}^j x_{i-2}^j+{\cal B}_{i}^j x_{i-1}^j+
{\cal C}_{i}^j x_{i}^j+ {\cal D}_{i}^j x_{i+1}^j+
{\cal E}_{i}^j x_{i+2}^j={\cal F}_i^j\
\]
for $i=1,...,n\,$ subject to periodic b.c.
$x_{-1}^j=x_{n-1}^j,\ \ x_0^j=x_n^j,\ \ x_{n+1}^j=x_1^j,\ \ x_{n+2}^j=x_2^j$
where
\begin{eqnarray*}
&&{\cal A}_i^j = \frac 1{r_{i}^j q_{i-1}^j r_{i-1}^j},\ 
{\cal E}_i^j = \frac 1{r_{i+1}^j q_{i+1}^j r_{i+2}^j},\  
\nonumber\\
&&{\cal B}_i^j = -\left(\frac 1{r_i^j q_{i-1}^j r_{i-1}^j}
+\frac 1{(r_i^j)^2 q_{i-1}^j}+\frac 1{(r_i^j)^2 q_i^j}
+\frac 1{r_i^j q_i^j r_{i+1}^j}\right)+
\nonumber\\
&&\hskip 1truecm + \ \frac{\phi(k_{i}^{j})+\phi(k_{i-1}^{j})}{2 r_i^j}+
\frac{\alpha_{i}^j}{2}- \frac34 \frac{ (k_{i+1}^{j})^2-(k_i^{j})^2}{q_i^j}
\nonumber\\
&&{\cal D}_i^j = -\left(\frac 1{r_i^j q_{i}^j r_{i+1}^j}
+\frac 1{(r_{i+1}^j)^2 q_{i}^j}+\frac 1{(r_{i+1}^j)^2 q_{i+1}^j}
+\frac 1{r_{i+1}^j q_{i+1}^j r_{i+2}^j}\right)+
\nonumber\\
&&\hskip 1truecm + \ \frac{\phi(k_{i+1}^{j})+\phi(k_i^{j})}{2 r_{i+1}^j}-
\frac{\alpha_{i}^j}{2}+\frac 34 \frac{(k_{i+1}^{j})^2-(k_i^{j})^2}{q_i^j}\,,
\nonumber
\\
&&
{\cal C}_i^j = \frac{q_i^j}{\tau}
-({\cal A}_i^j+{\cal B}_i^j+{\cal D}_i^j+{\cal E}_i^j),\   
{\cal F}_i^j = \frac {q_i^j}{\tau}x_i^{j-1},\nonumber\\
\end{eqnarray*}
where $\phi(k)=k^2 - \frac{b(k)}{k}$.
The initial quantities for the algorithm are computed from a discrete
representation of the initial curve $x_0$, for details see \cite{MS_CVS}.
Every pentadiagonal system is solved by mean of Gauss-Seidel iterations.
Next we present results of numerical simulations for the curve evolution driven 
by (\ref{geomrov}). In our experiments evolving curves are 
represented by $n=100$ grid points and we use discrete time step $\tau=0.001$. 
First we numerically compute time evolution of the initial 
ellipse with the halfaxes ratio 2:1 for the case  $b(k)=0$ (surface diffusion flow). We consider the time 
interval $[0,2]$. The evolution of curves without considering tangential
redistribution indicates accumulation of some curve representing grid points
and a poor resolution in other parts of the asymptotic shape, see also Figure
\ref{fig1}, a). In the case of asymptotically uniform tangential 
redistribution, 
we can see a uniform 
discrete resolution of the asymptotic shape, see Figure \ref{fig1}, b). 
Next we present evolution of an nontrivial initial curve driven by
(\ref{geomrov}) with $b(k)=0$. We show evolution of a highly nonconvex initial curve 
(see Figure \ref{fig1}, c) with  asymptotically uniform  redistribution 
($\omega=1$).
Since elastic curve dynamics is very fast in case of highly varying
curvature along the curve we have chosen smaller time step $\tau=10^{-6}$.
In Figure~\ref{astroid-test} we present evolution 
of an initial asteroid driven by (\ref{geomrov}) 
with $b(k)=-\frac12 k^3$ (the Willmore flow). 
A solution computed by the direct Lagrangian method 
is depicted by cross marks whereas solid curves correspond 
to the solution computed by the level set method approach. 
For details we refer the reader to \cite{BMOS}.

\section*{Acknowledgments}
This research was supported by grants: VEGA 1/3321/06, APVV-RPEU-0004-07 
(K.Mikula) and APVV-0247-06 (D.\v{S}ev\v{c}ovi\v{c}).

\bibliographystyle{ws-procs9x6}
\bibliography{mikula_karol}

\end{document}